\documentclass[12pt]{amsart}

\usepackage[utf8]{inputenc}			
\usepackage[T1]{fontenc}			
\usepackage[british]{babel} 		

\usepackage{bm,		
amsmath,			
amssymb,			
amsthm,				
latexsym,			
mathrsfs,			
mathtools,			
hyperref,			
geometry,			
}

\newtheorem{theorem}{Theorem}
\newtheorem{proposition}[theorem]{Proposition}
\newtheorem{cor}[theorem]{Corollary}
\newtheorem{lmm}[theorem]{Lemma}

\theoremstyle{definition}
\newtheorem{definition}[theorem]{Definition}
\newtheorem{remark}[theorem]{Remark}
\newcommand{\NN}{\mathbb N}
\newcommand{\sth}{\mathrel{:}}

\begin{document}

\title{Rough Approximate Subgroups}          
\author{Arturo Rodr\'\i guez Fanlo and Frank O. Wagner}
\address{Hebrew University of Jerusalem; Einstein Institute of Mathematics; Jerusalem, 9190401, Israel}
\address{Universit\'e Lyon 1; CNRS; Institut Camille Jordan UMR 5208, 21 avenue Claude Bernard, 69622 Villeurbanne Cedex, France}
\email{arturo.rodriguez@mail.huji.ac.il}
\email{wagner@math.univ-lyon1.fr}

\begin{abstract}
Given a rough definably amenable rough approximate subgroup $A$ of a group in some first-order structure, there is a type-definable subgroup $H$ normalised by $A$ and contained in $A^4$ of bounded index in $\langle A\rangle$.
\end{abstract}
\thanks{First author supported by the Israel Academy of Sciences and Humanities \& Council for Higher Education Excellence Fellowship Program for International Postdoctoral Researchers.}
\thanks{Second author supported by ANR-DFG AAPG2019 GeoMod}
\date{\today}

\maketitle

\section*{Introduction}
In \cite[Theorem 4.2]{H12}, Hrushovski found a fundamental connexion between approximate subgroups and Lie groups. This result, now known as the Lie Model Theorem, was the starting point for the complete classification of finite approximate subgroups by Breuillard, Green and Tao \cite{BGT12}. 

The proof of \cite[Theorem 4.2]{H12} first shows that, for any definable approximate subgroup $A$, under a certain definable amenability condition, there is a type-definable normal subgroup $S$ of bounded index of $\langle A\rangle$ contained in $A^4$. For this first fundamental step, Hrushovski proves a general Stabilizer Theorem \cite[Theorem 3.5]{H12}, working with an $S_1$ ideal of definable sets which, in addition, is invariant under translations and automorphisms. Then, in \cite[Theorem 4.2]{H12}, this general Stabilizer Theorem is applied for the $S_1$ ideal of sets of measure zero. However, in order for this ideal to become invariant under automorphisms, one needs to work with an expansion of the initial language.

In a joint work with Massicot, the second named author provided in \cite{MW15} an alternative proof of Hrushovski's Lie Model Theorem avoiding any expansion of the language. While Hrushovski's original proof only uses the ideal of measure zero subsets, the argument in \cite{MW15} fundamentally exploits the additivity of the measure.

In \cite{RF22}, the first named author adapted \cite[Theorem 3.5]{H12} and \cite[Theorem 4.2]{H12} for hyperdefinable sets, providing in particular a generalisation of the Lie Model Theorem for rough approximate subgroups. This result was later used in \cite{HR22} for the particular case of metric groups. The argument in \cite{RF22} uses again a general version of the Stabilizer Theorem (now for hyperdefinable sets), and so the applications in \cite{HR22} require again to work with an expansion of the initial language.

The aim of this paper is to adapt the alternative technique developed in \cite{MW15} for the case studied in \cite{HR22}. This is not a straightforward task --- while the ideal used for the original Lie Model Theorem is directly given by a measure, the ideal in \cite{HR22} is not associated with a measure at all. Instead, the ideal in \cite{HR22} is built as the inf-limit of the zero ideals of a sequence of subadditive set functions. Hence, in \cite{HR22}, we lack the additivity property that is fundamentally used in \cite{MW15}. Fortunately, a careful analysis of the properties of these subadditive set functions shows that, although they are not completely additive, they still keep some rough additivity that suffices to conclude the proof. \medskip

Our main result is the following:

\begin{theorem} \label{t:main} Assume $A$ is $T$-rough definably amenable with respect to $(\mu_i)_{i\in\NN}$ and suppose that for every $m\in\NN$ there is $i_m\in\NN$ such that $\mu_i(A^m)<\infty$ for all $i>i_m$. Then, there is a type-definable normal subgroup $N$ of $\langle AT_\omega\rangle$ of bounded index such that $T_\omega\subseteq N\subseteq A^4T_\omega$.
\end{theorem}

As a corollary, we get the following improved version of \cite[Theorem 20]{HR22}:

\begin{cor} \label{c:metric lie model} Let $(G_m,A_m,r_{i,m})_{i\leq m\in\NN}$ be a sequence such that, for some $\ell\in\NN$ and $K=(K_i\sth i\in\NN)$ sequence of natural numbers,
\begin{enumerate}
\item $G_m$ is a metric group,
\item $A_m$ is an $(\ell,r_{0,m})$-Lipschitz symmetric subset, 
\item $(r_{i,m})_{i\leq m}$ is a sequence of positive reals with $\max\{l,2\}\cdot r_{i,m}\leq r_{i-1,m}$ and
\[\mathrm{N}_{r_{i,m}}(A^4_m)\leq K_i\cdot \mathrm{N}_{r_{i,m}}(A_m)<\infty .\]
\end{enumerate}
Let $(G,A,r)$ be a non-principal ultraproduct in the group language with a predicate for $A$ and predicates $T_i$ for each ball of radius $r_i$. Let $T_\omega=\bigcap_{i\in\NN} T_i$. Then:\begin{enumerate}
\item There is a type-definable normal subgroup $N$ of $\langle AT_{\omega}\rangle$ of bounded index such that $T_\omega\subseteq N\subseteq A^4T_\omega$.\smallskip
\item If $K$ is constant, $A^2$ is a $T_\omega$-rough $K$-approximate subgroup. In particular, $AT_\omega$ has a connected Lie model $\pi:\ H\rightarrow L$, i.e.\ a surjective group homomorphism to a connected Lie group $L$ defined on a subgroup $H\leq \langle AT_\omega\rangle$ with kernel $T_\omega\trianglelefteq \ker(\pi)\subseteq A^8T_\omega$ such that
\begin{itemize}
\item $H\cap A^4T_\omega$ and $A^2T_\omega$ are commensurable,
\item $H\cap A^8T_\omega$ generates $H$, 
\item $\pi(H\cap A^4)$ is a compact neighbourhood of the identity in $L$, and
\item $\pi$ is continuous and proper from the logic topology using countably many parameters.
\end{itemize}\end{enumerate}
\end{cor}

\textbf{Notations:} Through this paper, fix a many sorted first order language $\mathtt{L}$, a $\kappa$-saturated $\mathtt{L}$-structure $\mathfrak{M}$ with $\kappa\geq |\mathtt{L}|$, a definable group $G$ and a definable symmetric subset $A$ of $G$. We always use product notation and, unless otherwise stated, we consider the group acting on itself on the left. A subset $X$ of $G$ is called symmetric if $1\in X=X^{-1}$. For subsets $X$ and $Y$ of $G$, we write $XY\coloneqq \{xy\sth x\in X,\, y\in Y\}$ for the set of pairwise products and call it the $Y$-thickening of $X$. For an integer $n>0$, abbreviate $X^n\coloneqq XX^{n-1}$ and $X^{-n}\coloneqq (X^{-1})^n$, and put $X^0\coloneqq\{1\}$. Write $Y^{X}\coloneqq\{x^{-1}yx\sth x\in X,\, y\in Y\}$, and say that $X$ normalises $Y$ if $Y^X\subseteq Y$. 

\section{Rough Measures}
Through this section, let $T\subseteq G$ be a symmetric subset and $\mathcal{A}$ some subalgebra of the boolean algebra of definable subsets of $G$ containing $A$ and closed under $T$-thickenings.

\begin{definition} A \emph{$T$-rough measure} on $\mathcal{A}$ is an increasing subadditive set function $\mu:\, \mathcal{A}\to\mathbb R_{\ge0}\cup\{\infty\}$ which is additive modulo $T$, i.e.\ for $X,Y\in\mathcal{A}$ with $X\cap YT=\emptyset$ we have $\mu(X\cup Y)=\mu(X)+\mu(Y)$. For a type-definable set $X$ we put $\mu(X)\coloneqq\inf\{\mu(\bar X)\sth\bar X\supseteq X\mathrm{\ definable}\}$. We say that $\mu$ is \emph{normalised} at $A$ if $\mu(A)=1$.\end{definition}

\begin{lmm}\label{l:union} Let $\mu$ be a $T$-rough measure on $\mathcal{A}$. Then, for any $S_0,\ldots,S_{n-1}\in\mathcal{A}$, we have
\[\mu\left(\textstyle{\bigcup_{i<n}}S_i\right)\ge\sum_{i<n}\mu(S_i)-\sum_{\mathclap{i<j<n}}\min\{\mu(S_i\cap S_jT),\mu(S_iT\cap S_j)\}.\]\end{lmm}
\begin{proof} For $n<2$, there is nothing to prove. For $n=2$, by subadditivity,
\[\mu(S_0)\le\mu(S_0\setminus S_1T)+\mu(S_0\cap S_1T).\]
Hence, by additivity modulo $T$, we get
\[\begin{aligned}\mu(S_0\cup S_1)&\ge\mu((S_0\setminus S_1T)\cup S_1)=\mu(S_0\setminus S_1T)+\mu( S_1)\\
&\ge\mu(S_0)+\mu(S_1)-\mu(S_0\cap S_1T).\end{aligned}\]
By symmetry, $\mu(S_0\cup S_1)\ge\mu(S_0)+\mu(S_1)-\mu(S_0T\cap S_1)$, and the result follows.

Suppose the assertion is true for $n$, and consider $S_0,\ldots,S_n\in\mathcal{A}$. Put
\[I=\{i<n\sth\mu(S_i\cap S_nT)\le\mu(S_iT\cap S_n)\},\]
and set $X=\bigcup_{i\in I}S_i$ and $X'=\bigcup_{i\in n\setminus I}S_i$. Suppose
$x\in X'\cap(S_n\setminus X'T)T$; say $x=yt$ with $y\in S_n\setminus X'T$ and $t\in T$. Then, $y=xt^{-1}\in X'T$, a contradiction. Therefore, $(X\cup X')\cap(S_n\setminus X'T)T\subseteq X\cap S_nT$. Thus,
\[\begin{aligned}\mu\left(\textstyle{\bigcup_{i\le n}}S_i\right)&\ge\mu(X\cup X'\cup(S_n\setminus X'T))\\
&\ge\mu(X\cup X')+\mu(S_n\setminus X'T)-\mu((X\cup X')\cap(S_n\setminus X'T)T)\\
&\ge\mu(X\cup X')+\mu(S_n)-\mu(S_n\cap X'T)-\mu(X\cap S_nT)\\
&\ge\sum_{i\le n}\mu(S_i)-\sum_{\mathclap{i<j<n}}\min\{\mu(S_i\cap S_jT),\mu(S_iT\cap S_j)\}\\
&\phantom{\geq}-\sum_{\mathclap{i\in n\setminus I}}\mu(S_iT\cap S_n)-\sum_{i\in I}\mu(S_i\cap S_nT)\\
&=\sum_{i\le n}\mu(S_i)-\sum_{\mathclap{i<j\le n}}\min\{\mu(S_i\cap S_jT),\mu(S_iT\cap S_j)\}.\end{aligned}\]
\end{proof}
\begin{remark} The above lemma also holds for type-definable sets, i.e.\ bounded intersections of definable sets in $\mathcal{A}$, by taking limits.\end{remark}

\begin{definition} Let $X$ be a definable subset of $G$. A definable symmetric subset $Y$ of $G$ is \emph{$t$-thick in $X$} if, for any $t+1$ elements $g_0,\ldots,g_t$ of $X$, there are $0\le i\neq j\le t$ such that $g_i^{-1}g_j\in Y$. A type-definable symmetric subset of $G$ is \emph{thick in $X$} if every definable symmetric superset is $t$-thick in $X$ for some $t\in\NN$.\end{definition}
Note that, since $Y$ is symmetric, $g_i^{-1}g_j\in Y$ if and only if $g_j^{-1}g_i\in Y$.
\begin{lmm} \label{l:thick covering} If $Y$ is $t$-thick in $X$, then $t$ left translates of $Y$ cover $X$. Conversely, if $t$ left translates of $Y$ cover $X$, then $Y^{-1}Y$ is $t$-thick in $X$.
\end{lmm}
\begin{proof} Take $Z\subseteq X$ maximal such that $z_1^{-1}z_2\notin Y$ for $z_1,z_2\in Z$ with $z_1\neq z_2$. Since $Y$ is $t$-thick, $Z$ exists and $|Z|\leq t$. Since $Z$ is maximal, for any $g\in X$, there is $z\in Z$ such that $z^{-1}g\in Y$; in other words, $X\subseteq ZY$.

Conversely, if $X\subseteq ZY$ with $|Z|\leq t$, by the pigeonhole principle, for any $g_0,\ldots,g_t$ in $X$, there are $0\le i\neq j\le t$ and $z\in Z$ such that $g_i,g_j\in zY$, so $g_i^{-1}g_j\in Y^{-1}Y$. Since $g_0,\ldots,g_t$ are arbitrary, we conclude that $Y^{-1}Y$ is $t$-thick
\end{proof}

\begin{lmm} \label{l:thick filter} Let $X$ be a definable subset of $G$.  An arbitrary intersection of type-definable subsets of $G$ which are thick in $X$ is again thick in $X$.\end{lmm}
\begin{proof} Let us see first that a finite intersection of definable thick subsets in $X$ is thick in $X$. If $Y\subseteq G$ is definable and $t$-thick in $X$ and $Y'\subseteq G$ is definable and $t'$-thick in $X$, then the intersection $Y\cap Y'$ is $R(2,t+1,t'+1)$-thick in $X$, where $R({\scriptscriptstyle{\bullet}},{\scriptscriptstyle{\bullet}},{\scriptscriptstyle{\bullet}})$ is the Ramsey number with respect to a $3$-colouring. Indeed, $Y\cap Y'$ is obviously definable and symmetric. Now, given any $Z$ subset of $X$ of size $R(2,t+1,t'+1)$, colour the edges of the complete graph on $Z$ by
\begin{itemize}
\item $\{z,z'\}$ is yellow if $z^{-1}z'\in Y\cap Y'$,
\item $\{z,z'\}$ is blue if $z^{-1}z'\notin Y$,
\item $\{z,z'\}$ is red if $z^{-1}z'\in Y\setminus Y'$.\end{itemize}
Since $Y$ is $t$-thick, $Z$ cannot contain a complete blue subgraph of size $t+1$ and, as $Y'$ is $t'$-thick, $Z$ cannot contain a complete red subgraph of size $t'+1$ either. Thus, it must contain a yellow edge, i.e.\ there are $z,z'\in Z$ with $z\neq z'$ such that $z^{-1}z'\in Y\cap Y'$. As $Z$ is arbitrary, we conclude that $Y\cap Y'$ is $R(2,t+1,t'+1)$-thick in $X$. By induction, from the two case, we conclude that any finite intersection of definable thick subsets in $X$ is again thick in $X$. 

Now, let $\bigcap_I Y_i$ be an intersection of type-definable subsets of $G$ thick in $X$ and $\bar Y$ be a definable symmetric superset of $\bigcap_I Y_i$. By compactness, there are a finite $I_0\subseteq I$ and definable symmetric supersets $\bar Y_i$ of $Y_i$ for each $i\in I_0$ such that $\bigcap_{I_0}\bar Y_i\subseteq\bar Y$. As $\bigcap_{I_0}\bar Y_i$ is thick in $X$, so is $\bar Y$; as $\bar Y$ is arbitrary, $\bigcap_I Y_i$ is thick in $X$.\end{proof}

\begin{lmm}\label{l:cover} Let $\mu$ be a $T$-rough measure on $\mathcal{A}$ normalised at $A$. Suppose that $A^{2m}\in\mathcal{A}$ and $\mu(A^{2m})<\infty$. Let $t>0$ be an integer and $B\in\mathcal{A}$ with $B\subseteq A^m$ and $\mu(B)\ge 2\mu(A^{2m})/t$. Put
\[\mathrm{S}(B)=\left\{g\in A^{2m}\sth\min\{\mu(B\cap gBT),\mu(gB\cap BT)\}\ge\frac{2\mu(A^{2m})}{t^2}\right\}.\]
Then, $\mathrm{S}(B)$ is $t$-thick in $A^m$.
\end{lmm}
Note that $\mathrm{S}(B)$ is symmetric, since $A$ is symmetric and $\mu$ is left invariant. Also, note that any $g\in G$ with $B\cap gBT$ or $gB\cap BT$ non-empty must be in $BTB^{-1}$, so, in particular, $\mathrm{S}(B)\subseteq BTB^{-1}$.
\begin{proof} Suppose, aiming for a contradiction, that the conclusion is false, and consider a counter-example $(g_i\sth i\le t)$. Then,
\[\min\{\mu(g_i B\cap g_jBT),\mu(g_i BT\cap g_jB)\}<\frac{2\mu(A^{2m})}{t^2},\]
for all $0\le i<j\le t$. By Lemma \ref{l:union}, 
\[\begin{aligned}\mu(A^{2m})&\ge\mu\Big(\textstyle{\bigcup_{i\le t}}\, g_iB\Big)\\
&\ge\sum_{i\le t}\mu(g_iB)-\sum_{\mathclap{0\le i<j\le t}}\min\{\mu(g_i B\cap g_jBT),\mu(g_i BT\cap g_jB)\}\\
&>(t+1)\cdot\frac{2\mu(A^{2m})}{t}-\frac{(t+1)t}{2}\cdot\frac{2\mu(A^{2m})}{t^2}\\
&=\mu(A^{2m})\cdot\frac{t+1}{t}>\mu(A^{2m}),\end{aligned}\]
a contradiction.
\end{proof}

\section{Rough Definable Amenability}
Throughout this section, let $T=(T_i\sth i\in \NN)$ be a sequence of definable symmetric subsets of $G$ such that $T_{i+1}^2\subseteq T_i$ and $T_{i+1}^A\subseteq T_i$ for all $i\in\NN$. Write $T_{\omega}=\bigcap_{i\in\NN} T_i$; note that $T_\omega$ is a type-definable subgroup of $G$ normalised by $A$. 

\begin{definition} We say that $A$ is \emph{$T$-rough definably amenable} if for all $i\in\NN$ there is a $T_i$-rough measure $\mu_i$ on the boolean algebra of definable subsets of $G$ invariant under left translations and normalised at $A$.\end{definition}
\begin{remark} We do not really need $\mu_i$ to be defined on all definable subsets of $G$, but just on a boolean subalgebra $\mathcal{A}$ which contains $A$ and is closed under left translations, $T_i$-thickenings for $i\in\NN$, and setwise products (i.e.\ $XY\in \mathcal{A}$ for any $X,Y\in \mathcal{A}$).\end{remark}

\begin{lmm}\label{l:sanders} Let $P$ be a partially ordered set and $f_0,\ldots,f_{r-1}:\, P\rightarrow \mathbb{R}_{>0}$ be increasing functions. Suppose there is $a>0$ such that $f_k(x)\ge a$ for all $x\in P$ and each $k<r$. Then, for any $\varepsilon_0,\ldots,\varepsilon_{r-1}>0$, there is $x\in P$ such that $f_k(y)>(1-\varepsilon_k)\cdot f_k(x)$ for any $y\in P$ with $y<x$ and each $k<r$.
\end{lmm}
\begin{proof} Aiming for a contradiction, suppose otherwise. Then, we can recursively construct a decreasing sequence $(x_n)_{n\in\NN}$ in $P$ such that for any $n\in\NN$ there is $i<r$ with $f_i(x_{n+1})\leq (1-\varepsilon_i)\cdot f_i(x_n)$. Since $r\in\NN$ is finite, by the pigeonhole principle, there is $i<r$ and $N\subseteq\NN$ infinite such that $f_i(x_{n+1})\leq (1-\varepsilon_i)\cdot f_i(x_n)$ for any $n\in N$. Since $f_i$ is monotonous, we get that $f_i(x_m)\leq (1-\varepsilon_i)\cdot f_i(x_n)$ for any $n\in N$ and $m>n$. Let $n_0\coloneqq \min N$ and $b\coloneqq f_i(x_{n_0})$. For any $n\in N$, put $\alpha\coloneqq |\{m\in N\sth m<n\}|$. Then,
\[a\leq f_i(x_n)\leq (1-\varepsilon_i)^\alpha\cdot f_i(x_{n_0})=(1-\varepsilon_i)^\alpha\cdot b,\]
getting a contradiction when $\alpha>\ln(a/b)/\ln(1-\varepsilon_i)$. \end{proof}

\begin{proposition}\label{p:square} Let $A\subseteq G$ be $T$-rough definably amenable, and $m$ and $r$ integers with $\mu_i(A^{2m(r+1)})<\infty$ for each $i\in\NN$. Suppose $B\subseteq A^m$ is type-definable and thick in $A$. Then, there is a type-definable subset $S\subseteq B^2T_\omega\cap A^{2m}$ thick in $A^m$ such that $S^{r}\subseteq B^4T_\omega$.\end{proposition}
\begin{proof}
For all $i,t\in\NN$, define recursively on $j\in\NN$ a family $\mathscr{T}^t_{i,j}$ of definable subsets of $A^m$. For $j=0$, let  $X\in\mathscr{T}^t_{i,0}$ if $X\neq\emptyset$. Given $\mathscr{T}^t_{i,j}$, let $X\in\mathscr{T}^t_{i,j+1}$ if 
\[S^t_{i,j}(X)\coloneqq\left\{g\in A^{2m}\sth X\cap gXT_i,X\cap g^{-1}XT_i\in \mathscr{T}^{t^2}_{i,j}\right\}\]
is $t$-thick in $A^m$. Note that $S^t_{i,j}(X)\subseteq XT_iX^{-1}$ is symmetric.

By induction on $n$, we easily see that these families are definable, i.e.\ $\{b\sth X(b)\in\mathscr{T}^t_{i,j}\}$ is definable for any definable family $X(b)$ of subsets of $A^m$. Also, we easily see by induction that $Y\in \mathscr{T}^t_{i,j}$ implies $X\in\mathscr{T}^s_{i,\ell}$ and $S^t_{i,j}(Y)\subseteq S^s_{i,\ell}(X)$ whenever $j\ge \ell$, $t\le s$ and $Y\subseteq X$. 

Next, we claim that, if $X\subseteq A^m$ is definable with $\mu_i(X)\ge 2\mu_i(A^{2m})/t$, then $X\in\mathscr{T}^t_{i,j}$ for all $j\in\NN$. This is trivial for $j=0$. Assume it holds for $j$; by Lemma \ref{l:cover}, we have that 
\[\mathrm{S}_i(X)\coloneqq \{g\in A^{2m}\sth \min\{\mu_i(X\cap gXT_i),\mu_i(X\cap g^{-1}XT_i)\}\geq 2\mu_i(A^{2m})/t^2\}\]
is $t$-thick in $A^m$ and, by induction hypothesis, $\mathrm{S}_i(X)\subseteq S^{t^2}_{i,j}(X)$, concluding that $X\in\mathscr{T}^t_{i,j+1}$.\smallskip

Let $\mathcal{B}$ be the set of definable symmetric supersets of $B$ contained in $A^m$. Fix $\bar{B}\in\mathcal{B}$ and $i\in \NN$. By Lemma \ref{l:thick covering}, finitely many left translates of $\bar B$ cover $A$, so we have $\mu_i(\bar B)>0$. Let $K_i=\mu_i(A^{2m(r+1)})$ and $\varepsilon_i=\mu_i(\bar B)/2rK_i$. Finally, set $t_{i,n}\coloneqq \lceil 2K_i/\mu_i(\bar B)\rceil^{2^n}$, so $\bar B\in\mathscr{T}^{t_{i,0}}_{i,j}$ for all $j\in\NN$.

We define recursively sequences $(X_{i,n}\sth n\in\NN)$ and $(S_{i,n}\sth n\in\NN)$ such that $(X_{i,n}\sth n\in\NN)$ is decreasing, $S_{i,n}$ is thick in $A^m$, $S_{i,n}\subseteq X_{i,n}T_iX_{i,n}\cap A^{2m}$ is symmetric and $X_{i,n}\in\mathscr{T}^{t_{i,n}}_{i,j}$ for all $j\in\NN$. 

For $n=0$, set $X_{i,0}=\bar B$. Given $X_{i,n}$, set $S_{i,n}=\bigcap_{j\in\NN}S^{t_{i,n}}_{i,j}(X_{i,n})$; since $X_{i,n}\in\mathscr{T}^{t_{i,n}}_{i,j}$ for all $j\in\NN$, it follows that $S_{i,n}$ is thick in $A^m$ by Lemma \ref{l:thick filter}. Also, note that $S_{i,j}(X_{i,n})\subseteq X_{i,n}T_iX_{i,n}^{-1}\cap A^{2m}$ is symmetric. Given $X_{i,n}$ and $S_{i,n}$, we take $X_{i,n+1}=X_{i,n}\cap gX_{i,n}T_i$, where $g\in S_{i,n}$ satisfies for some $k<r$ that
\[\mu_i(X_{i,n+1}T_i^k\bar BT_i^k\cap A^{2m(r+1)})\le(1-\varepsilon_i)\mu_{i}(X_{i,n}T_i^k\bar BT_i^k\cap A^{2m(r+1)}),\]
or, if there is no such $g\in S_{i,n}$, we take $X_{i,n+1}=X_{i,n}$. Since $X_{i,n}\in\mathscr{T}^{t_{i,n}}_{i,j}$ and $g\in S^{t_{i,n}}_{i,j}(X_{i,n})$ for all $j\in\NN$, we conclude that $X_{i,n+1}\in\mathscr{T}^{t_{i,n+1}}_{i,j}$ for all $j\in\NN$.
\smallskip

Since $1\in T_i$ and $X_{i,n}\subseteq \bar B\subseteq A^m$ for every $n\in\NN$, we have by monotonicity and left invariance that 
\[0<\mu_i(\bar B)\leq \mu_i(X_{i,n+1}T_i^k\bar B T_i^k\cap A^{2m(r+1)})\leq \mu_i(A^{2m(r+1)})< \infty.\] 
Applying Lemma \ref{l:sanders} to the functions $f_k:\ X\mapsto \mu_i(XT_i^k \bar B T_i^k\cap A^{2m(r+1)})$, $k<r$, on the set $\{X_{i,n}\sth n\in\NN\}$ linearly ordered by inclusion, we conclude that there is $n_i\in \NN$ such that 
\[\mu_i(X_{i,n}T_i^k \bar{B}T_i^k\cap A^{2m(r+1)})>(1-\varepsilon_i)\cdot \mu_i(X_{i,n_i}T_i^k \bar{B}T_i^k \cap A^{2m(r+1)})\]
for all $n>n_i$ and each $k<r$. By construction of the sequence $(X_{i,n})_{n\in\NN}$, we see that $X_{i,n_i}=X_{i,n_i+1}\eqqcolon X_i$, and
\[\mu_i((X_i\cap gX_iT_i)T_i^k \bar{B}T_i^k\cap A^{2m(r+1)})>(1-\varepsilon_i)\cdot \mu_i(X_iT_i^k \bar{B}T_i^k \cap A^{2m(r+1)})\]
for all $g\in S_i:=S_{i,n_i}$ and each $k<r$. Thus
\[\begin{aligned}\mu_i(X_iT_i^k \bar{B}T_i^k\cap gX_iT_i^{k+1}\bar{B}T_i^k\cap A^{2m(r+1)})
&\geq\mu_i((X_i\cap gX_iT_i)T_i^k \bar{B}T_i^k\cap A^{2m(r+1)})\\
&>(1-\varepsilon_i)\cdot \mu_i(X_iT_i^k \bar{B}T_i^k\cap A^{2m(r+1)})
\end{aligned}\]
for all $g\in S_i$ and each $k<r$. Hence,
\[\begin{aligned}\mu_i((&X_iT_i^k \bar{B}T_i^k\setminus gX_iT_i^{k+1} \bar{B} T_i^{k+1})\cap A^{2m(r+1)})\\
&\leq \mu_i(X_iT_i^k \bar{B} T_i^k\cap A^{2m(r+1)})-\mu_i(X_iT_i^k \bar{B} T_i^k\cap gX_iT_i^{k+1} \bar{B} T_i^k\cap A^{2m(r+1)})\\
&<\big(1-(1-\varepsilon_i)\big)\cdot\mu_i(X_iT_i^k \bar{B} T_i^k\cap A^{2m(r+1)})\\
&\leq\varepsilon_i\,K_i.\end{aligned}\]

Take $g_0,\ldots,g_ {r-1}\in S_i$ arbitrary and write $g_{<k}\coloneqq g_0\cdots g_{k-1}$ for each $k<r$. It follows that
\[\begin{aligned}\mu_i((X_i \bar{B} &\setminus g_{<r}X_iT_i^{r} \bar{B} T_i^{r})\cap A^{2m})\\
&\leq\mu_i\Big(\bigcup_{k<r}g_{<k}(X_iT_i^k \bar BT_i^k\setminus g_{k} X_iT_i^{k+1} \bar{B} T_i^{k+1})\cap A^{2m}\Big)\\
&\leq\sum_{k<r}\mu_i(g_{<k}(X_iT_i^k\bar{B}T_i^k\setminus g_{k} X_iT_i^{k+1} \bar{B}T_i^{k+1})\cap A^{2m})\\
&\leq\sum_{k<r}\mu_i((X_iT_i^k \bar{B}T_i^k\setminus g_{k} X_iT_i^{k+1} \bar{B} T_i^{k+1})\cap A^{2m(r+1)})\\
&<r\,\varepsilon_i\,K_i<\mu_i(\bar{B})\le\mu_i(X_i \bar{B}\cap A^{2m}).\end{aligned}\]
In particular, $X_i\bar{B}\cap g_{<r} X_iT_i^{r}\bar{B}T_i^{r}\neq\emptyset$, so $g_{<r}\in X_i \bar{B}T_i^r\bar{B}T_i^rX_i\subseteq \bar{B}^2(T_i^r\bar{B})^2$. \smallskip

Now, put $S(\bar{B})\coloneqq \bigcap_{i\in\NN}S_i$ and $S=\bigcap_{\bar{B}\in\mathcal{B}}S(\bar{B})$. Then, for any $g_0,\ldots,g_{r-1}\in S$,  we conclude that
\[g_{<r}\in\bigcap_{\bar{B}\in\mathcal{B}}\bigcap_{i\in\NN}\bar{B}^2(T_i^r\bar{B})^2.\]
By saturation, $\bigcap_{\bar{B}\in\mathcal{B}}\bigcap_{i\in\NN} \bar{B}^2(T_i^r\bar{B})^2=B^2(T^r_\omega B)^2=B^4T_\omega$ --- recall that $T_\omega$ is a subgroup normalised by $A$, whence also normalised by $B$. Hence, $S^{r}\subseteq B^4T_\omega\subseteq A^{4m}T_\omega$. Also note that $S\subseteq B^2T_\omega\cap A^{2m}$. \end{proof}

We can now prove Theorem \ref{t:main}:
\begin{theorem}\label{t:approximate} Assume $A$ is $T$-rough definably amenable and suppose that for every $m\in\NN$ there is $i_m\in\NN$ such that $\mu_i(A^m)<\infty$ for all $i>i_m$. Then, there is a type-definable normal subgroup $N$ of $\langle AT_\omega\rangle$ of bounded index such that $T_\omega\subseteq N\subseteq A^4T_\omega$.
\end{theorem}
\begin{proof} We iterate Proposition \ref{p:square}. Put $A_0=A$ and, given $A_n\subseteq A^{2^n}$, we apply Proposition \ref{p:square} with $2^n$ in place of $m$, $8$ in place of $r$, the subsequences $(\mu_i\sth i>i_{9\cdot 2^{n+1}})$ and $(T_i\sth i>i_{9\cdot 2^{n+1}})$ in place of $(\mu_i)_{i\in\NN}$ and $(T_i)_{i\in\NN}$ respectively, $A_n$ in place of $B$ and set $A_{n+1}$ in place of $S$. Note that $A_{n+1}\subseteq A_n^2\subseteq A^{2^{n+1}}$. 

Clearly, $H=\bigcap_{n\in\NN}A_n^4T_\omega$ is a type-definable subgroup contained in $A^4T_\omega$. Furthermore, $H$ is thick in $A^m$ for each $m\in\NN$. Thus, the index of $H$ in $\langle AT_\omega\rangle$ is bounded by Lemma \ref{l:thick covering}. Now, as $H$ has bounded index in $\langle AT_\omega\rangle$, the intersection of all $\langle AT_\omega\rangle$-conjugates of $H$ is a bounded intersection, whence a type-definable normal subgroup $N$ in $\langle AT_\omega\rangle$ of bounded index.
\end{proof}

\begin{definition} Let $K>0$ be an integer and $Z\subseteq G$ a symmetric subset. We say that $A\subseteq G$ is a {\em $Z$-rough $K$-approximate
subgroup} if there is a finite set $E$ of size $K$ such that $A^2\subseteq EAZ$.\end{definition}
\begin{remark} The notion of rough approximate subgroups has been studied by several authors before (e.g.\ \cite{T14}, \cite{GL20}). Typically, rough approximate subgroups are only defined in the context of metric groups where $Z$ is a ball at the identity of positive radius. Our definition is more flexible as it does not requires a metric. Nevertheless, note that, when $Z$ is symmetric, both definitions coincide as one can (a posteriori) pick the discrete metric given by $\mathrm{d}(x,y)=\min\{n\sth x^{-1}y\in Z^n\}$.
\end{remark}
\begin{lmm} \label{l:approximate} Assume $A$ is $T$-rough definably amenable, and suppose $\mu_{i+3}(A^4)\leq K_i$ for some $i\in\NN$. Then, $A^2$ is a $T_i$-rough $K_i^2$-approximate subgroup.\end{lmm}
\begin{proof} This is just a rough version of Ruzsa's Covering Lemma. Suppose $(g_k\sth k\le K_i)$ are in $A^3$ with $g_kA\cap g_{k'}AT_{i+3}=\emptyset$ for $k\neq k'$. Hence, 
\[K_i\ge \mu_{i+3}(A^4)\ge (K_i+1)\mu_{i+3}(A)=K_i+1,\]
a contradiction. Thus, there is $E\subseteq A^3$ of size at most $K_i$ such that, for any $g\in A^3$, there is $e\in E$ with $gA\cap eAT_{i+3}\neq\emptyset$, whence $A^3\subseteq EAT_{i+3}A\subseteq EA^2T_{i+2}$.

It follows that 
\[(A^2)^2=A^3\cdot A\subseteq EA^2T_{i+2}\cdot A=EA^3T_{i+1}\subseteq E^2A^2T_{i+2}T_{i+1}\subseteq E^2A^2T_i,\]
which yields that $A^2$ is a $T_i$-rough $K_i^2$-approximate subgroup.\end{proof}

\begin{cor} \label{c:approximate} Assume $A\subseteq G$ is $T$-rough definably amenable, and suppose that $\mu_i(A^4T_\omega)<\infty$ for all $i\in\NN$. Then, $\mu_i(A^m)<\infty$ for all $i,m\in\NN$. In particular, by Theorem \ref{t:approximate}, there is a type-definable normal subgroup of bounded index of $\langle AT_\omega\rangle$ contained in $A^4T_\omega$.
\end{cor}
\begin{proof} By Lemma \ref{l:approximate}, $A^2$ is a $T_j$-rough approximate subgroup for all $j\in\NN$. Thus, for any $m,j\in\NN$, there is $E_{m,j}$ finite such that $A^m\subseteq E_{m,j}A^2T_j$. Since $\mu_i(A^2T_\omega)<\infty$, there is $j_0\in\NN$ such that $\mu_i(A^2T_{j_0})<\infty$, so $\mu_i(A^m)\leq |E_{m,j_0}|\cdot \mu_i(A^2T_{j_0})<\infty$. \end{proof}

\section{Applications to Metric Groups}
In this section we explain how to apply Theorem \ref{t:approximate} to get Corollary \ref{c:metric lie model} (i.e.\ \cite[Theorem 20]{HR22}). Recall that a \emph{metric group} is a group $G$ together with a metric $\mathrm{d}$ invariant under left translations. We denote by $\bar{\mathbb{D}}_r$ the closed ball of radius $r$ at the identity. For $\ell\in \NN_{>0}$ and $r\in\mathbb{R}_{>0}$, a subset $A\subseteq G$ is \emph{$(\ell,r)$-Lipschitz} if every right translation by an element of $A$ is $\ell$-Lipschitz in $\bar{\mathbb{D}}_r$, i.e.\ $\mathrm{d}(xa,ya)\leq \ell\cdot\mathrm{d}(x,y)$ for any $a\in A$ and $x,y\in \bar{\mathbb{D}}_r$. For a subset $X\subseteq G$, we define $\mathrm{N}_r(X)\coloneqq \max\{|Z|\sth Z\subseteq X$ and $\mathrm{d}(z,z')>r$ for all $z,z'\in Z$ with $z\neq z'\}$. 
\begin{remark} Note that $\mathrm{N}_{r}$ is a $\bar{\mathbb{D}}_r$-rough measure invariant under left translations. Indeed, as the metric is invariant under left translations, so is $\mathrm{N}_{r}$. If $Z\subseteq X\cup Y$ is $r$-separated, $|Z|\leq |Z\cap X|+|Z\cap Y|$ where $Z\cap X\subseteq X$ and $Z\cap Y\subseteq Y$ are $r$-separated, so $\mathrm{N}_r$ is subadditive. Finally, if $X\cap Y\bar{\mathbb{D}}_{r}=\emptyset$, then $X$ and $Y$ are $r$-separated. As the union of two $r$-separated subsets which are $r$-separated is $r$-separated, $\mathrm{N}_{r}(X\cup Y)=\mathrm{N}_{r}(X)+ \mathrm{N}_{r}(Y)$, so $\mathrm{N}_r$ is $\bar{\mathbb{D}}_r$-rough additive. 
\end{remark}
Fix natural numbers $\ell$ and $(K_i\sth i\in\NN)$. Let $(G_m,A_m,r_{i,m})_{i\leq m\in\NN}$ be such that
\begin{enumerate}
\item $G_m$ is a metric group,
\item $A_m$ is an $(\ell,r_{0,m})$-Lipschitz symmetric subset,
\item $r_m=(r_{i,m})_{i\leq m}$ is a finite sequence of positive real numbers with $\max\{\ell,2\}\cdot r_{i,m}\leq r_{i-1,m}$ and
\[\mathrm{N}_{r_{i,m}}(A^4_m)\leq K_i\cdot \mathrm{N}_{r_{i,m}}(A_m)<\infty .\]
\end{enumerate}
Consider the language of groups enriched with a predicate for $A_m$ and predicates for $\bar{\mathbb{D}}_{r_{i,m}}$ for each $i\in\NN$ --- where $r_{i,m}=r_{m,m}$ for $i\geq m$; we consider $(G_m,X_m,r_{i,m})_{i\leq m\in\NN}$ as a sequence of structures in this language. Pick a non-principal ultrafilter $\mathfrak{u}$ of $\NN$ and let $(G,A,r)$ be the corresponding ultraproduct --- recall that this is $\aleph_1$-saturated \cite[Exercise 5.2.3]{TZ12}. Set $T_i\coloneqq \bar{\mathbb{D}}_{r_i}$ and $T_\omega\coloneqq\bigcap T_i$. 

\begin{lmm} $A$ is $T$-rough definably amenable with respect to $\mu_i(Y)\coloneqq \frac{\mathrm{N}_{r_i}(Y)}{\mathrm{N}_{r_i}(A)}$ and $\mu_i(A^4T_\omega)<\infty$ for any $i\in\NN_{>0}$. 
\end{lmm}
\begin{proof} As $A$ is $(\ell,r_0)$-Lipschitz and $\ell r_{i+1}\leq r_i$, we get $T_{i+1}^A\subseteq T_i$. Since $2r_{i+1}\leq r_i$, we have $T_{i+1}T_{i+1}\subseteq T_i$. As $\mathrm{N}_{r_i}$ is a  $T_i$-rough measure on the boolean algebra of definable subsets invariant under left translations, we conclude that $A$ is $T$-rough definably amenable with respect to the normalisation $\mu_i$. Finally, note that $\mu_i(A^4T_\omega)\leq \mu_i(A^4T_{i+2})\leq \mu_{i-1}(A^4)\leq K_{i-1}$ for any $i\in\NN_{>0}$. 
\end{proof} 

By Corollary \ref{c:approximate}, there is a type-definable normal subgroup $N$ of $\langle AT_\omega\rangle$ of bounded index such that $T_\omega\subseteq N\subseteq A^4T_\omega$; in other words, we get Corollary \ref{c:metric lie model}(1).\smallskip

When $(K_i\sth i\in\NN)$ is a constant $K$, $A^2$ is a $T_\omega$-rough $K$-approximate subgroup by Lemma \ref{l:approximate} and saturation. Thus, by \cite[Theorem 3.2]{RF22}, we get Corollary \ref{c:metric lie model}(2).

\end{document}